\documentclass [11pt]{article}
\usepackage{amssymb}
\usepackage{amsmath}

\newcommand{\N}{\mathbb{N}}

\newcommand{\Z}{\mathbb{Z}}

\newtheorem{definition}{{\bf Definition}}[section]
\newtheorem{theorem}[definition]{{\bf Theorem}}

\newtheorem{claim}[definition]{\noindent {\bf Claim}}

\newtheorem{lemma}[definition]{\noindent {\bf Lemma}}
\def\endproof{\hfill {\kern 6pt\penalty 500
\raise -0pt\hbox{\vrule \vbox to5pt {\hrule width 5pt
\vfill\hrule}\vrule}}}
\def\centerpicture #1 by #2 (#3){\leavevmode
        \vbox to #2{
        \hrule width #1 height 0pt depth 0pt
        \vfill
        \special{pictfile #3}}}
\makeindex
\begin{document}
\title{Small clones and the projection property}

\author{Maurice Pouzet\footnote{Research done under the auspices of Intas programme 03-51-4110 "Universal algebra and lattice theory"}\\
\small  Math\'ematiques, ICJ\\
\small Universit\'e Claude-Bernard, Lyon1\\
\small  Domaine de Gerland -b\^at. Recherche [B],\\
\small 50 avenue Tony-Garnier, \\
\small F$69365$ Lyon, France\\ 
\small {\rm e-mail: pouzet@univ-lyon1.fr} \\
\small Fax 33 4 37 28 74 80 \and Ivo G.Rosenberg\\
\small D\'ept. de math\'ematiques  et de statistique,\\
\small Universit\'e de Montr\'eal,\\
\small CP 6128 succ. centre-Ville,\\
\small Montr\'eal, QC, Canada,H3C3J7\\
\small {\rm e-mail: rosenb@dms.umontreal.ca}\\  
\small Fax 1-514 341 5700}  

\date{\today}

\maketitle


\begin{abstract} In $1986$, the second author classified the minimal clones on a finite universe into five types. We extend this classification to infinite universes and to multiclones.  We show that every non-trivial clone contains a "small" clone of one of the five types. From it
we deduce, in part,  an earlier result, namely that if $\mathcal C$  is a clone on a universe $A$ with at least two elements, that contains all constant
operations, then all binary idempotent operations are projections and some $m$-ary idempotent
operation is not a projection  some $m\geq 3$ if and only if there is a Boolean group
$G$ on $A$ for which $\mathcal C$ is the set of all operations $f(x_1,\dots, x_n)$ of the form
$a+\sum_{i\in I}x_i$ for $a\in A$ and
$I\subseteq \{1,\dots, n\}$.  
\end{abstract}
Keywords: Clones and  multiclones.

\section{Definitions, notations, results} Let $A$ be a fixed universe of cardinality at least $2$. Denote by $\mathfrak P(A)$ the family of subsets of $A$. Denote by $\N$ the set of
non-negative integers and  by $\N^*$ the set $\N\setminus \{0\}$.  For $n\in \N^*$ a map $f: A^n \rightarrow \mathfrak P(A)$ is
an
$n$-{\it ary multioperation} on $A$. Call $f$ an \emph{hyperoperation} if the empty set is not in $im f$, the image of $f$. Special binary hyperoperations, called \emph{hypergroups}, where introduced in 1934 \cite{marty}  and there is a sizeable literature on them, see \cite{corsini91}, \cite{corsini03}. If $im f\subseteq \{\{b\}: b\in A\}\cup \{\emptyset\}$, the map $f$ is called a \emph{partial operation} on $A$. If we identify each singleton value $\{b\}$ of $f$ with the element $b$ and declare that $f$ is not defined  at all $a\in A^n$ with $f(a)=\emptyset$ then, clearly, $f$ becomes a partial operation in the usual sense. If $f$ is such that $im f\subseteq \{\{b\}: b\in A\}$ the same identification of all singleton values  yields an operation on $A$ in the usual sense. Denote by $\mathcal M^{(n)}$, $\mathcal H^{(n)}$, $\mathcal P^{(n)}$ and $\mathcal O^{(n)}$ the set of
$n$-ary multioperations, hyperoperations, partial operations and operations  on $A$.  Set $\mathcal M:=\bigcup \{\mathcal M^{(n)}: n\in \N^* \}$ and define similarly $\mathcal H, \mathcal P$ and $\mathcal O$. For a subset $\mathcal Z$ of $\mathcal M, \mathcal H, \mathcal P$ and $ \mathcal O$,  call the pair $<A; \mathcal Z>$ a (nonindexed) \emph{multialgebra, hyperalgebra, partial algebra} and \emph{algebra} on $A$. 

For $\mathcal C\subseteq \mathcal M$ and $n\in \N^*$ set $\mathcal
C^{(n)}:=\mathcal C\cap  \mathcal M^{(n)}$. In the sequel
$\approx$ denotes an identity on $A$ (i.e. the two sides are equal for all values of the variables  in $A$) and we write $:\approx $ for
a defining identity. For $f\in \mathcal M^{(n)}$ and $\pi$ a permutation of the set $\{1,\dots,n\}$, the operation $f_{\pi}$ defined by
$f_{\pi}(x_{1},\dots,x_{n}):\approx f(\pi(x_{1}),\dots, \pi(x_{n})) $ is an  {\it isomer} of $f$. A  multioperation
$f\in
\mathcal M$ is\ {\it idempotent} provided
$f(x,\dots,x)\approx \{x\}$. We denote by $\mathcal I$ the set of all idempotent multioperations on $A$. For
$n=1,2,3$ the multioperations from
$\mathcal M^{(n)}$ are called, respectively, {\it unary}, {\it binary} and {\it ternary}.    A ternary multioperation $f$ is a {\it majority}
multioperation    if
$f(x,x,y)\approx f(x,y,x)\approx f(y,x,x)\approx \{x\}$;  if, on the opposite, $f(x,x,y)\approx f(x,y,x)\approx f(y,x,x)\approx \{y\}$  the
multioperation $f$  is a {\it minority} multioperation.  If $f(x,x,y)\approx y \approx f(y,x,x)$ the
multioperation $f$ is a {\it
Mal'tsev multioperation}  and, if in addition $f(x,y, x)\approx \{x\}$, this is a {\it Pixley multioperation}.  For
$i,  n\in \N^*$ with
$i\leq n$, define the
$i$-{\it th
$n$-ary projection}
$e^{n}_{i}$ by setting $e^{n}_{i}(x_{1},\dots,x_{n}):\approx x_{i}$. Set 
$\mathcal Q:= \{e^{n}_{i}: i,n \in \N^*\}$. For $n\geq 3$ and $1\leq i \leq n$,
call $f\in \mathcal M^{(n)}$ a {\it semiprojection on its $i$-th coordinate} if $f(a_{1},\dots,a_{n})=\{a_{i}\}$ whenever
$a_{1},\dots,a_{n}\in A$ are not pairwise distinct.

On the set $\mathcal M$ we define a countable set $\{\pi_{ij}: i,j\in \N^*\}$ of partial operations. For $i,j\in \N^*$, the map $\pi_{ij}: \mathcal M^{(i)}\times (\mathcal M^{(j)})^{i}\rightarrow \mathcal M^{(i)}$ is defined as follows. Let $f\in   \mathcal M^{(i)}$ and $(g_1,\dots, g_i)\in (\mathcal M{(j)})^{i}$. For every $a:=(a_1, \dots, a_j)\in A^{j}$ set:
\begin{equation}\label{multioperation}
\pi_{ij}(f, g_1, \dots, g_i)(a):=\bigcup \{f(u_1,\dots, u_i): u_1\in g_1(a), \dots,  u_i\in g_i(a)\}
\end{equation} 

Notice that (\ref{multioperation}) makes sense since $g_1, \dots, g_i)(a)$ and $f(u_1,\dots, u_i)$ are subsets of $A$. Also, if $f,g_1, \dots, g_i$ are operations, the right-hand side of (\ref{multioperation}) is the standard $f(g_1(a),\dots, g_i(a))$.

 For $\mathcal Z\subseteq \mathcal M$ denote by $[\mathcal Z]$ the least member of $\mathcal M$ containing $\mathcal Z\cup \mathcal Q$ and closed under all $\pi_{ij}$, with $i,j\in \N^{*}$. The subsets of $\mathcal M$ of the form $[\mathcal Z]$ for some $\mathcal Z \subseteq \mathcal M$ are called \emph{multiclones}. If $\mathcal Z\subseteq \mathcal X$ for $\mathcal X\in \{\mathcal H, \mathcal P, \mathcal O\}$, then $[\mathcal Z]$ is said to be, respectively, a \emph{hyperclone}, \emph{partial clone} and \emph{clone}.  For  $\mathcal Z\subseteq \mathcal O$ the clone  $[\mathcal Z]$ is the set of term operations of the algebra $<A, \mathcal Z>$; this fact can be extended to $\mathcal Z\subseteq \mathcal X\in \{\mathcal M, \mathcal P, \mathcal H\}$.

A multiclone is said to be {\it minimal} if  $\mathcal Q$ is its only proper submulticlone.
An operation $f$ is said to be {\it minimal} if the multiclone $[f]$ generated by $f$ is minimal and $f$ is of minimum arity among
the multioperations in $[f]\backslash \mathcal Q$. For an example, let $f_n\in \mathcal M^{(n)}$ with $f_n(a):=\emptyset$ for all $a\in A^n$. Then $\mathcal Q\cup \{f_n: n\in \N^{*}\}$ is a minimal multiclone  and $f_1$ is minimal.  It is known that on  a finite universe, every clone distinct from $\mathcal Q$ 
contains a minimal clone and, as shown by the second author \cite {rosfive},  the minimal operations fall into  five 
distinct types.  Despite the fact that on an infinite universe a clone distinct  from $\mathcal Q$ does not necessarily contain  
 a minimal clone, it turns out that the main feature of the second author's result  is preserved: clones, and particularly minimal ones,
can be classified by means of the five types of operations, \'A. Szendrei \cite{szendrei}. It was mentioned by J.~Pantovi\'c and G.Vojvodic \cite{pantovic} that on  a finite universe, every hyperclone   distinct from $\mathcal Q$ 
contains a minimal hyperclone of one of the five types. Here we show first that the classification into five types extends to multiclones on an arbitrary universe.

In order to present our result, we recall that a {\it Boolean group} is a
2-elementary  group $<A; +,0>$, that is a group with neutral element $0$ satisfying $a+a=0$ for all $a\in A$. It is well known that a
Boolean group is necessarily abelian; in fact such a group on  $A$ finite exists exactly if $\vert
A\vert$ is a power of
$2$ and in that case the group is isomorphic to  a power of
$\Z/2\Z$.

\begin{theorem} \label {smallclone}Let $\mathcal C$ be a multiclone on $A$. If $\mathcal C\not=\mathcal Q$,
then  $\mathcal C\setminus\mathcal Q$ contains either:\\
$1)$ a unary multioperation;\\
$2)$ a binary  idempotent multioperation;\\
$3)$ a majority multioperation;\\
$4)$ a semiprojection;\\
$5)$ the term operation $x+y+z$ of a Boolean group $<A; +,0>$.\\
 \end{theorem}

Next, we apply our result to the projection property, a property we introduced in  \cite {pou-ros-sto} for structures of various sorts, like posets, graphs, metric spaces.
A structure 
$R$ is 
$n$-{\it projective} if the only idempotent morphisms from  its $n$-power $R^{n}$ into $R$  are the projections. We looked
at the relationship between these properties for various values of $n$.  One of our   results ( Theorem 1.1 of
\cite {pou-ros-sto}, see also \cite{pou-ros-sto2}) can be deduced, in part, from Theorem \ref{smallclone}

\begin {theorem}\label{projection}
The following are equivalent for a   clone $\mathcal C$ on a universe $A$ with at least two
elements. \\
$(i)$  $\mathcal C$    contains all constant operations,  all its binary idempotent operations are projections while some 
$n$-ary idempotent operation is not a projection;\\ $(ii)$ there is a  Boolean group $G$ on $A$
for which $\mathcal C$ is the set $F_{G}$ of all operations  of the form
$f(x_1,\dots, x_n)\approx a+\sum_{i\in I}x_i$ for $a\in A$, $I\subseteq \{1,\dots, n\}$ and $n\in \N^{*}$.  
\end{theorem}
A result similar to Theorem \ref{projection} was obtained independently in \cite {dav-all} (Lemma 2.6 ).

\section {Proof of Theorem \ref{smallclone}}
Our proof is an adaptation of the proof
from
\cite {rosfive} (see also \cite {pou-ros-sto}).\\

Let $n$ be the least positive integer such that $\mathcal C^{(n)}\not= \mathcal Q^{(n)}$. Notice that for $n>1$, $f\in
\mathcal C^{(n)}$ and $1\leq i<j\leq n$, the multioperation $g(x_{1},\dots, x_{n-1}):\approx
f(x_{1},\dots,x_{i}, \dots, x_{j-1}, x_{i},x_{j}, \dots,x_{n-1})$ belongs to   $\mathcal
C^{(n-1)}$ and hence is a projection. If
$n=1, 2$  then any 
$f
\in
\mathcal C^{(n)}\setminus  \mathcal Q^{(n)}$ satisfies  $1)$ or   $2)$ of Theorem
\ref{smallclone}.  If $n
\geq 4$ then, according to the following lemma, any 
$f
\in
\mathcal C^{(n)}\setminus  \mathcal Q^{(n)}$ satisfies $4)$.
  
\begin {lemma} Let $n, n\geq 4$, and  $f\in \mathcal C^{(n)}$; suppose  that every multioperation
obtained from
$f$ by identifying two variables is a projection, then $f$ is a semiprojection.
\end{lemma}
\begin{proof}  We assume $\vert f(a)\vert=1$ for all $a:=(a_1, \dots, a_n)\in A^n$ such that $a_1, \dots, a_n$ are not pairwise distinct. We can then apply the well-known Swierczkowski Lemma\cite {swier}. Indeed, as it turns out,  the fact that we may have $\vert f(a)\vert\not =1$, for some $a:=(a_1, \dots, a_n)\in A^n$ with $a_1, \dots a_n$ pairwise distinct,  is irrelevant.
\end{proof}

Thus, we may  suppose that $n=3$.  In this case, the following lemma asserts that we may find $f$ in case $3$, $4$ or $5$, which
proves Theorem \ref{smallclone}.
 
 \begin{lemma} \label {main} Suppose $\mathcal C^{(2)}=\mathcal Q^{(2)}$ and $\mathcal C^{(3)}\setminus
\mathcal Q^{(3)}\not=\emptyset$. If
$\mathcal C^{(3)}\setminus \mathcal Q^{(3)}$  contains no semiprojection and no majority multioperation then $\mathcal
C^{(3)}\setminus \mathcal Q^{(3)}=\{m\}$, where $m$ is a totally symmetric minority multioperation. If, moreover,
$\mathcal C^{(4)}\setminus
\mathcal Q^{(4)}$  contains no semiprojection then {m} is the term operation
$x+y+z$ of a Boolean group $<A; +,0>$.
\end{lemma}

\begin{proof} 

Let $f\in \mathcal C^{(3)}\setminus \mathcal Q^{(3)}$. As  every multioperation
obtained from
$f$ by identifying two variables is a projection, there are $a,b,c\in \{1,2\}$ such that 
\begin{equation} \label{eq 1}
f(x_{1},x_{1},
x_{2})\approx  x_{a},  \;\; \; f(x_{1},x_{2}, x_{1})\approx  x_{b},\; \; \; f(x_{1},x_{2}, x_{2})\approx  x_{c}
\end{equation}
We denote by $\chi_{f}$ the ordered triple $abc$ and  abreviate $(1)$ by: 
\begin{equation} \label{eq 2}
\overline f (112)=a, \; \; \;\overline f (121)=b,\;\;\;\overline f (122)=c.
\end{equation} 
Thus,  if  $\chi_{f}\in \{111, 122, 212\}$ then $f$ is a semiprojection; if $\chi_{f}=112$  then $f$ is
a majority multioperation; if $\chi_{f}=221$ then $f$ is a minority multioperation, and 
if  $\chi_{f} =211$ then $f$ is a Pixley multioperation.\\

\begin{claim} \label{Claim 1.} {\it If $\chi_{f}\in \{121,222\}$ then $\chi_{h}=211$ for some $h\in
\mathcal C^{3}$}.
\end{claim}
 \noindent{\bf Proof of Claim \ref{Claim 1.}.}  If $\chi_{f}=121$,   set $h(x_{1},x_{2},
x_{3}):\approx  f(x_{1},x_{3}, x_{2})$. 
 Then $$\overline h(112)=\overline f(121)=2, \;\; \overline h(121)=\overline f(112)=1, \;\;  \overline h(122)=\overline f(122)=1$$
proving $\chi_{h}=211$.   If
$\chi_{f}=222$,  set $h(x_{1},x_{2}, x_{3}):\approx f(x_{2},x_{1}, x_{3})$.  Then  $\overline
h(112)=\overline h(211)= \overline h(212)=2$ and so $ \overline h(112) =2$, $\overline h(121)=1$ and $\overline h(121)=1$,
proving $\chi_{h}=211$ and the claim.\endproof\\

As it is well known (see e.g. Theorem 9.3.2 p.201 \cite{deneckewismath}) a clone
$\mathcal C$  contains a Pixley operation if and only if it contains a majority and a Mal'tsev
operation; in fact, as shown in the proof of Lemma 2.7 \cite{rosfive} p. 413, this amounts to
the fact that
$\mathcal C$ contains a majority and a minority operation. This extends to multioperations.  For reader's convenience, we reprove what we need.\\ 

\begin{claim}\label{Claim 2.} {If  $\chi_{h}=211$ then $\mathcal C$ contains a majority multioperation.}
\end{claim} 

\noindent{\bf Proof of  Claim \ref{Claim 2.}.} Set $m(x_{1}, x_{2}, x_{3}):\approx h(x_{1}, h(x_{1},  x_{2}, x_{3}),
x_{3})$.
In view of $\chi_{h}=211$, we  get:
$$\overline m(112)=\overline h(1\overline h (112)2)=\overline h (112)=1,$$
$$\overline m(121)=\overline h(1\overline h (121)1)=\overline h  (111)=1, $$
$$\overline m(122)=\overline h(1\overline h (122)2)=\overline h ( 112)=1\overline
h (111)=1.$$ proving $\chi_{g}=221$  and the claim.
 \endproof

Supposing that $\mathcal C^{(2)}=\mathcal Q^{(2)}$ and that $\mathcal C^{(3)}\setminus
\mathcal Q^{(3)}$ is non-empty and contains no semiprojection and no majority multioperation, it follows from Claim \ref{Claim 1.} and Claim \ref{Claim 2.} 
that $\mathcal C^{(3)}\setminus
\mathcal Q^{(3)}$ contains only minority multioperations.

\begin{claim} \label{Claim 3.} {\it Let $f_{1}, f_{2}\in \mathcal C^{(3)}$ be  minority multioperations. Then }
\begin{equation}
\label{eq 3} f_{1}(f_{2}(x_{1},x_{2},x_{3}),x_{2},x_{3})\approx x_{1}
\end{equation}
\begin{equation}
\label{eq 4} f_{1}(x_1, f_{2}(x_{1},x_{2},x_{3}),x_{3})\approx x_{2}
\end{equation}
\begin{equation}
\label{eq 5} f_{1}(x_1, x_2, f_{2}(x_{1},x_{2},x_{3}))\approx x_{3}
\end{equation}
\end{claim}
\noindent {\bf Proof of Claim \ref{Claim 3.}.} Denote by $s(x_{1},x_{2},x_{3})$ the left-hand side of $(\ref {eq 3})$. Then we have successively:
$$\overline s(112)=\overline f_{1}(\overline f_{2}(112)12)=\overline f_{1}(212)=1,$$ 
$$\overline s(121)=\overline f_{1}(\overline f_{2}(121)21)=\overline f_{1}(221)=1,$$ 
$$\overline s(122)=\overline f_{1}(\overline f_{2}(122)22)=\overline f_{1}(122)=1.$$
Hence, $s$ is a semiprojection.  Since $s\in \mathcal C^{(3)}$ and  $\mathcal C^{(3)}\setminus \mathcal Q^{(3)}$ contains
no semiprojection, $s$ is a projection, that is  $s=e^{3}_{1}$ proving $(\ref {eq 3})$. Identities $(\ref {eq 4})$ and $(\ref {eq 5})$ follow from $(\ref {eq 3})$
applied to the isomers $f_1(x_2, x_1, x_3)$ and $f_1(x_2, x_3, x_1)$ of $f$. This proves  the claim.\endproof

\begin{claim}\label{Claim 4.} {\it Let $f\in \mathcal C$ be a minority multioperation and $a,b,c\in A$.  Then $f(a,b,c)=\{a\}$  if and only if $b=c$.}\end{claim}
\noindent {\bf Proof of Claim \ref{Claim 4.}.} If $b=c$ then, since $f$ is a minority, $f(a,b,c)=f(a,b,b)=\{a\}$. Conversely, suppose $f(a,b,c)=\{a\}$. From
(\ref{eq 4}) of Claim \ref{Claim 3.}  applied to $f_1=f_2:=f$, we have $ f(x_1, f(x_{1},x_{2},x_{3}),x_{3})\approx x_{2}$, in particular $ f(a,
f(a,b,c),c)= \{b\}$. With $f(a,b,c)=\{a\}$ , this gives  $f(a, a, c)=\{b\}$.  Since $f$ is a minority, we get $b=c$, as claimed. \endproof

\begin{claim}\label {Claim 5.} {\it $\mathcal
C$ contains only one minority multioperation.}\end{claim}

\noindent{\bf Proof of Claim \ref{Claim 5.}.} Let $f_{1},f_{2}\in \mathcal C$ be minority multioperations. Set 
$$h(x_1, x_2, x_3):\approx f_{1}(x_1, f_1 (x_1, x_2, x_3),  f_2(x_1, x_2, x_3)).$$ 

Proceeding as in the above proof of (\ref {eq 3}) we obtain that 
$h$ is a semiprojection on the first coordinate, hence the first projection, that is
  $$f_1(x_1, f_1 (x_1, x_2, x_3),  f_2(x_1, x_2,
x_3))\approx x_1.$$ According to Claim \ref{Claim 4.}
$$f_1(x_1,x_2, x_3)
\approx f_{2} (x_1, x_2, x_3).$$

Since all isomers of a minority multioperation are minority multioperations, it follows from  the uniqueness of a minority multioperation $g\in
\mathcal C$ that $g$ is totally symmetric, that is invariant under all permutations of variables. This completes the proof of the
first part of Lemma \ref {main}.

For simplicity, we write $(x_{1}x_{2}x_{3})$ for $g(x_{1},x_{2},x_{3})$.\\
From now, we suppose that  $\mathcal C^{(4)}\setminus \mathcal Q^{(4)}$ contains no semiprojection. For fixed $a,b\in A$, define $\varphi_{ab}:AÊ\rightarrow \mathfrak{P}(A)$ by $\varphi_{ab}(x):\approx g(x,a,b)$.
\begin{claim}\label {Claim 5a.}{\it For all $a,b\in A$ the map $\varphi:= \varphi_{ab}$ is a permutation of $A$ of order at most $2$ (i.e. an involution) and consequently $g\in \mathcal O^{3}$ (i.e. an operation). } \end{claim}

\noindent{\bf Proof of Claim \ref{Claim 5a.}.} From (\ref{eq 3}) for $f_1=f_2=g$ we have $$\varphi^2(x)=g(g(x,a,b),a,b)\approx \{x\}.$$

We show that $\varphi\in \mathcal O^{(1)}$, that is $\varphi $ is a selfmap of $A$. Indeed, let $c\in A$ and
$u\in \varphi(c)$ be arbitrary. Then $\varphi(u)\in  \varphi^2(c)=\{c\}$ and thus $\{u\}= \varphi^2(u)=\varphi(c)$ proving that $\vert \varphi(c)\vert =1$.  It follows that $\varphi$ is a permutation of $A$ such that $\varphi^2=e^1_1$.\endproof

\begin{claim}\label {Claim 6.}{\it Fix $0\in A$ and put $x+y:\approx (xy0)$. Then  $<A;+, 0>$ is a Boolean
group and $(xyz)\approx x+y+z$. } \end{claim}
{\bf Proof of Claim \ref{Claim 6.}. } From the total symmetry of $( \;)$ and from $x 00\approx x$ we obtain that $<A; +, 0>$ is a
commutative groupoid with the neutral element $0$. Next, $a+b=0$ if and only if $a=b$. Indeed, if $a+b=0$, then
$(ab0)=0$, hence $a=b$ by Claim \ref{Claim 4.}. Conversely, we have $a+a=(aa0)=0$.

We prove that $(xyz)\approx
x+(y+z)$, that is, in view of the observation just above,  $(xyz)+(x+(y+z))\approx 0$. Using the definition of our groupoid
operation, this means 
$((xyz)(x(yz0)0)0)\approx 0$. It suffices then to prove that the  following identity holds 
\begin{equation}
\label{eq 6}
((xyz)(x(yzt)t)t)\approx t .
\end{equation}

Let $h$ be the quaternary term operation on $<A; (\;)>$ defined by :\\
\begin{equation}
\label{eq 7}
h(x_{1},x_{2},x_{3},x_{4}):\approx ((x_{1}x_{2}
x_{3})(x_{1}(x_{2}x_{3}x_{4}) x_{4})x_{4})).
\end{equation}
We first show  that $h$ is a semiprojection on its last variable. We abbreviate the right-hand side of (\ref{eq 7}) by $((123)
(1(234 )4)4)$. We consider  the six possibilities  of  identifying two variables. \\

$a)$. Set $x_{1}=x_{2}$. Using the fact that $(\;  )$ is a totally symmetric minority operation and  $(4)$ of Claim \ref{Claim 3.} (for
$f_{1}=f_{2}=(\; )$), we get  $$((113)( 1(134)4)4)=(3((314)14)4) =(3(1(134)4)4)=(334)=4.$$

$b)$. Set $x_{1}=x_{3}$. Due to the total symmetry of $(\;)$ this case reduces to the previous one.\\ 

$c)$. Set $x_{1}=x_{4}$.  Using the fact that $(\; )$ is a totally symmetric minority operation, we get
$$((423)
(4(234 )4)4)=((423)(234)4)=((234)(234)4)=4.$$

$d)$. Set $x_{2}=x_{3}$. Using the fact that $(\;)$ is a minority operation we get $$((122)(1(224 )4)4)=(1(144)4)=(114)=4.$$ 

$e)$. Set $x_{2}=x_{4}$. Using the fact that $(\;)$ is a totally symmetric minority operation we get 
$$((143)
(1(434 )4)4)=((143)(134)4)=((134)(134)4)=4.$$

$f)$. Set $x_{3}=x_{4}$. As above we get 
$$((124)
(1(244 )4)4)=((124)(124)4)=4.$$
Since $\mathcal C^{(4)}\setminus\mathcal Q^{(4)}$ contains no semiprojection,  $h$ is a projection, in fact $h=e^{4}_{4}$,
proving that  identity $(\ref{eq 6})$ holds. From this identity and the total symmetry of $(\;)$  we get 
$x+(y+z)\approx(zxy)\approx(xyz)\approx x+(y+z)$ proving that the binary operation $+$ is associative. This concludes the
proof of the claim. \endproof

The proof of Lemma \ref{main} is complete.
\end{proof}

\section{Proof of Theorem \ref{projection}} 
If $G$ is a $2$-elementary group and $\mathcal C= F_G$ (where $F_G$ was introduced in Theorem \ref{projection}) then, 
clearly,  $\mathcal C$ contains all constant maps, all its  binary idempotent members are projections and the term operation 
$x+y+z$ is idempotent and not a projection. 

\noindent Conversely,
suppose that
$\mathcal C$ contains all constant operations, that its  idempotent binary operations are the two projections and some $n$-ary
idempotent operation is not a projection  for a fixed $n\geq 3$. With the following lemma, we get that $\mathcal C\setminus
\mathcal Q$ contains neither a  semiprojection nor a
majority operation.
\begin{lemma}
Let  $\mathcal C$ be a clone on a universe $A$ with at least two
elements, that contains all constant operations and such that the binary idempotent operations are the two projections.  Then for
all $n\geq 2$ an
 operation
$g\in
\mathcal C^{(n)}$ is a projection if an only if some isomer $f$ of $g$ satisfies
\begin{equation}\label{eq 8}
f(y,y, x_{3},\dots, x_{n})\approx y.
\end{equation} 
\end{lemma}

\begin{proof}
 The proof is an adaptation of Lemma 2.4 of \cite{pou-ros-sto}. The necessity of (\ref{eq 8}) is obvious.
We prove the sufficiency of (\ref{eq 8}) by induction on $n$ ($n\geq 2$). For
$n=2$ if $f(y,y)\approx y$, $f$ is idempotent and $f$ is a projection by the hypothesis. Suppose $n\geq 3$ and every $h\in
\mathcal C^{(n-1)}$ satisfying (\ref{eq 8}) is $e_{i}^{n-1}$ for some $i\in \{1,2\}$. Since $\mathcal C$ contains all constant
operations, for every
$a\in A$,  it contains the $(n-1)$-ary operation
$f_{a}$  defined by 
$f_{a}(x_1,\dots, x_{n-1}):\approx f(x_1,\dots, x_{n-1},a)$. From the inductive hypothesis, $ f(x_1,\dots,x_{n-1},a)\approx
x_{i(a)}$ for some
$i(a)\in
\{1,2\}$. If
$i(a)=2$ for all
$a\in A$ then $f$ is the second projection. Suppose $i(b)=1$ for some $b\in A$. The $(n-1)$-ary operation $f'$ defined by
$f'(x_{1},\dots,x_{n-1}):\approx f(x_{1},\dots,x_{n-1}, x_{n-1})$  belongs to $\mathcal C$ and satisfies the same hypothesis as
$f$, hence it is either the first or the second projection. Since, from above, 
$f'(x_{1}, \dots, x_{n-2}, b)\approx f(x_{1}, \dots, x_{n-2}, b,b)\approx x_{1}$, clearly $f'$ is the first projection. Let $a\in A$
be arbitrary. From
$x_{1}\approx f(x_{1}, \dots, x_{n-2}, a, a )\approx x_{i(a)}$ we obtain $i(a)=1$ proving that $f$ is the first projection. This
proves the inductive step and the lemma.
\end{proof}\\
 
Let $\mathcal C':= \mathcal C\cap \mathcal I$ be the clone of the idempotent operations from $\mathcal C$. According to
Theorem
\ref{smallclone},
$\mathcal C'^{(3)}\setminus \mathcal Q^{(3)}=\{m_{G}\}$ where $m_{G}(x,y,z)\approx x+y+z$ for a Boolean group $G:= <A;
+,0>$. From \cite {pou-ros-sto} (the statement at the end  of page 173) follows  that
$\mathcal C=F_{G}$. With this the proof of Theorem \ref{projection} is complete.

\end{document}